\documentclass[11pt]{article}

\usepackage{amssymb,latexsym   }
\usepackage[centertags]{amsmath}

\numberwithin{equation}{section}

\newtheorem{thm}[equation]{Theorem}
\newtheorem{lemma}[equation]{Lemma}

\newtheorem{cor}[equation]{Corollary}

\def\done{{1\hskip-2.5pt{\rm l}}}

\begin{document}

\title{Zeroes of Gaussian analytic functions}

\author{Mikhail Sodin
\thanks{Supported by
the Israel Science Foundation of the Israel Academy of Sciences
and Humanities.}
\\ School of Mathematics
\\ Tel Aviv University
\\ Tel Aviv 69978 Israel
\\ sodin@post.tau.ac.il
}

\maketitle

\begin{abstract}
Geometrically, zeroes of a Gaussian analytic function are
intersection points of an analytic curve in a Hilbert space with a
randomly chosen hyperplane. Mathematical physics provides another
interpretation as a gas of interacting particles. In the last
decade, these interpretations influenced progress in understanding
statistical patterns in the zeroes of Gaussian analytic functions,
and led to the discovery of canonical models with
invariant zero distribution. We shall discuss some of
recent results in this area and mention several open questions.
\end{abstract}

\section*{Introduction}

A Gaussian analytic function is a linear combination
\[
f(z) = \sum_{k\ge 0} \zeta_k f_k(z)
\]
of analytic functions $f_k \colon G\to \mathbb C$ ($G \subseteq
\mathbb C$ is a domain),
\[
\sum_{k\ge 0} |f_k(z)|^2 < \infty \qquad\text{locally\ uniformly\
in}\ G,
\]
with independent standard complex Gaussian random coefficients
$\zeta_k$. The random zero set $\mathcal Z_f = f^{-1}(0)$ is the
theme of this talk.

Pioneering contributions in this area were made by Paley and
Wiener~\cite[Chapter~X]{PW}, Kac~\cite[Chapter~I]{Kac}, and
Rice~\cite{Rice}. Paley and Wiener constructed a large class of
Gaussian analytic functions in a strip with stationary
distribution with respect to the shifts, and computed the mean
number of zeroes. Their work was influenced by ergodic theory and
theory of almost-periodic functions. Kac was interested in the
mean number of real zeroes of polynomials with real coefficients.
Rice systematically treated both theoretical and applied aspects
of random noises in radio signals. The technique introduced by Kac
and Rice has been significant to radio engineers and physicists.

These studies were continued in various directions, notably by
Littlewood and Offord~\cite{LO}, Hammersley~\cite{Ham},
Offord~\cite{Offord, Offord1}, and Kahane
\cite[Chapter~13]{Kahane}. Hammersley looked at pure probabilistic
aspects of point processes generated by zeroes of random
polynomials. The other authors were motivated by the entire
functions and Nevanlinna theory. Introducing randomness, they
tried to single out `typical patterns' in the zero distribution in
various classes of entire and analytic functions.

In the 90s, the subject was revived by several groups of
researchers who came from different areas: Bogomolny, Bohigas and
Leboeuf; Shub and Smale; Edelman and Kostlan; Hannay; Bleher,
Shiffman and Zelditch; Nonnenmacher and Voros; by no means is this
list complete. They established new links to physics (Coulomb gas
of charged particles, random matrices, quantum chaos) and geometry
(analytic curves in projective Hilbert spaces), and drastically
changed the whole subject.

\medskip
We start this lecture with a quick review of basic results
pertaining to zero sets of {\em arbitrary} Gaussian analytic
functions (Section~\ref{sect1}). In Section~\ref{sect2}, we
introduce three canonical random zero process (on the Riemann
sphere, complex plane, and the unit disc), distinguished by
stationarity with respect to the isometries. In
Sections~\ref{sect3} and \ref{sect4}, we  consider in more detail
one of them, the canonical
random zero process in $\mathbb C$. The exposition in this part is
based on joint works with Tsirelson~\cite{ST}. In
Section~\ref{sect5}, we discuss an `exactly solvable case' of the
hyperbolic zero process recently discovered by Peres and
Vir\'ag~\cite{PV}.

\section{The random zero process $\mathcal Z_f$}\label{sect1}
Informally speaking, the random zero set $\mathcal Z_f=f^{-1}(0)$
is the intersection of an analytic curve $\mathfrak f\colon G \to
P(H)$ in a projective Hilbert space with a random hyperplane, the
analytic functions $f_k$ are homogeneous coordinates of the curve
$\mathfrak f$. Projective unitary transformations of the curve
$\mathfrak f$ do not change the random zero set $\mathcal Z_f$.
Hence the random set $\mathcal Z_f$ depends only on geometry of
the curve $\mathfrak f$. Its study can be interpreted as part of
the H.~Cartan -- Ahlfors -- H.~and J.~Weyl theory of analytic
curves independent of the dimension of the target space.
\footnote{An interesting attempt to build a `dimensionless theory'
of analytic curves was made by Favorov \cite{Favorov}. His
approach is based on the pluripotential theory in Banach spaces.}

The properties of the random process $\mathcal Z_f$ can be
expressed by its counting measure
\[
n_f = \sum_{a\colon f(a)=0} \delta_a\,,
\]
$\delta_a$ is a point measure at $a$. The measure $n_f$ is a
random, positive, locally finite measure on $G$. The classical
formula
\begin{equation}\label{eq.10}
n_f = \frac1{2\pi} \Delta \log|f|
\end{equation}
(the Laplacian is understood in the sense of distributions)
explicitly relates the random measure $n_f$ to the Gaussian
analytic process $f$. Proofs of most of the results presented
below start with this relation.

\subsection{The Edelman-Kostlan formula for the mean measure} The
first question about the random measure $n_f$ is to find its
average which is a non-negative measure in $G$.
\begin{thm}[Edelman-Kostlan \cite{EK}]\label{thm.EK}
\[
\mathbb E n_f = \frac1{2\pi} \Delta \log \|\mathfrak f\|\,,
\]
where
\[
\|\mathfrak f\|(z) := \sqrt{\sum_{k\ge 0} |f_k(z)|^2}\,.
\]
I.e., the mean measure $\mathbb En_f$ coincides with the Riesz
measure of the subharmonic function $\log \|\mathfrak f\|(z)$.
\end{thm}

The RHS of the Edelman-Kostlan formula is a pull-back of the
Fubini-Study area measure from the projective space $P(H)$ to $G$
by the curve $\mathfrak f$. Its density with respect to the
Euclidean area measure in $G$ equals
\begin{equation}\label{eq.FS}
\frac1{\pi} (\mathfrak f^\#)^2 := \frac1{\pi}
\frac{\sum_{i<k} |f'_i f_k - f_i f'_k|^2}{\|\mathfrak f\|^4}\,.
\end{equation}
The function $\mathfrak f^\#$ is a `Fubini-Study derivative' of
the curve $\mathfrak f$. The Edelman-Kostlan formula can be viewed
as a version of the classical Crofton formula from the integral
geometry. Its proof is a simple computation based on equation
\eqref{eq.10}:
\[
\mathbb E n_f = \frac1{2\pi} \Delta (\mathbb E \log |f|) =
\frac1{2\pi} \Delta \log \|f\|\,,
\]
since, for any complex Gaussian random variable $\zeta$,
\[
\mathbb E \log |\zeta| = \log \|\zeta \| + \text{const}\,.
\]

What about the higher `moments' of the random measure $n_f$? They
are expressed by the $k$-point correlation measures
\begin{equation}\label{eq.mu}
d\mu (z_1,\,...\, z_k) = \mathbb E \left( dn_f(z_1)\,...\,
dn_f(z_k)\right)
\end{equation}
on $\underbrace{G\times\,...\,G}_{k\ \text{times}}$. Hannay
\cite{Hannay96} derived explicit formulas for these measures which
generalize \eqref{eq.FS}. They involve determinants and permanents
of $k\times k$ matrices. In different contexts, the rigorous proof
of these formulas is given in \cite{BSZ} and \cite{PV}.

\subsection{Calabi's rigidity}

Surprisingly, the mean $\mathbb E n_f$ {\em determines}
the random zero set $\mathcal Z_f$. In geometry, the same
phenomenon was discovered by Calabi already in the beginning of the
50s.
\begin{thm}\label{thm.Calabi}
Let $f$ and $g$ be Gaussian analytic functions in a domain $G$,
and let $\mathbb E n_f = \mathbb E n_g$. Then the corresponding random
zero sets $\mathcal Z_f$ and $\mathcal Z_g$ have the same
distribution.
\end{thm}
This holds due to the underlying analyticity.
The idea is not difficult: let $K(z_1, z_2) = \mathbb E \left(
f(z_1)\overline{f(z_2)} \right)$ be the covariance of the process
$f$. By the Edelman-Kostlan formula, the mean measure $\mathbb E
n_f$ determines the function $z\mapsto \log K(z,z)$ up to a
harmonic summand. In turn, the diagonal $K(z,z)$ determines the
whole covariance kernel $K$ (due to analyticity of $K$ in $z_1$
and $\bar{z_2}$), and hence the distribution of the Gaussian
process $f$. The details and references are in \cite{Sodin}.

Here is Calabi's original formulation: {\em If two linearly
non-degenerate analytic curves $\mathfrak f\colon M\to P(H_1)$,
$\mathfrak g\colon M\to P(H_2)$ of a complex manifold $M$ induce
the same Riemannian metric on $M$ by pulling back the
corresponding Fubini-Study metrics, then the projective spaces
coincide $P(H_1) = P(H_2)$, and the curves are unitarily
equivalent.}

\subsection{Offord-type estimate}\label{sec.Offord}
\begin{thm}\label{thm.Offord}
Let $f$ be a Gaussian analytic function on
a plane domain $G$. Then for any test function $\phi\in C^2_0(G)$
with a compact support in $G$, and any $\lambda>0$
\begin{equation}\label{eq.20}
\mathbb P \left( \left| \int_G \phi (dn_f - \mathbb E(dn_f)\,
)\right| > \lambda \right) \le 3 \exp \left( -
\frac{2\pi\lambda}{\|\Delta \phi \|_1} \right)\,.
\end{equation}
\end{thm}
Here, $\|\,.\,\|_1$ is the $L^1$ norm with respect to the area
measure.

Here is an argument borrowed from  Offord \cite{Offord}. By
\eqref{eq.10} and Green's formula, we need to estimate the probability
\[
\mathbb P \left( \left| \int_G (\log|f| - \mathbb E \log |f| )
\Delta\phi\, dm \right| > 2\pi \lambda \right)\,,
\]
$m$ is the area measure. This reduces the proof to a simple fact
about concentration of $\log |\zeta|$, where $\zeta$ is a complex
Gaussian random variable. The details are in \cite{Sodin}.

The result can be extended in various directions. It persists for
zero sets of any random analytic process $f$ in $G$ with uniformly
bounded exponential moment: for some $c>0$,
\[
\sup_{z\in G} \mathbb E \left( e^{c |\log|f(z)|\,|}\right) <
\infty\,.
\]
Examples of such analytic processes are given in \cite{NSV}.
Instead of the $L^1$-norm of the Laplacian $\Delta\phi$, one can
fix any of the $L^q$-norms, $2<q\le \infty$, of the gradient
$\nabla\phi$. In Section~\ref{sect4}, we discuss a more
complicated `global version' of Theorem~\ref{thm.Offord}.

There is a price for such a level of generality: sometimes,
Offord's estimate does not give an optimal result. For example, it
does not yield sharp bounds for the `hole probability' (see
\eqref{eq.60} and \eqref{eq.70} below).

\section{Chaotic analytic zero points}\label{sect2}
In the beginning of the nineties, Bogomolny, Bohigas and Leboeuf;
Kostlan; and Shub and Smale introduced a remarkably unique class of
Gaussian
analytic functions with unitary invariance of
zero points. Following Hannay \cite{Hannay96}, we use the term
`chaotic analytic zero points' (CAZP, for short). We consider here
three CAZP models: the elliptic CAZP, the flat CAZP, and the
hyperbolic CAZP \footnote{The toric CAZP (which we do not discuss
here) was introduced by Nonnenmacher and Voros \cite{NV}.}. They
are the random zero set of a Gaussian analytic function
\begin{align}
f_L(z) &= \sum_{k=0}^L \zeta_k \sqrt{ \frac{L(L-1) \dots
(L-k+1)}{k!} } \, z^k && \text{(elliptic, $ L = 1, 2, \dots $)},
\label{elliptic} \displaybreak[0]
\\ f_L(z) &= \sum_{k=0}^\infty \zeta_k \sqrt{ \frac{L^{k}}{k!}} \,
z^k && \text{(flat, $ L>0 $)}, \label{flat} \displaybreak[0]
\\
f_L(z) &= \sum_{k=0}^\infty \zeta_k \sqrt{ \frac{L(L+1)
\dots(L+k-1)}{k!}}\, z^k && \text{(hyperbolic, $ L>0$)}.
\label{hyperbolic}
\end{align}
The analytic function (\ref{elliptic}) is a polynomial of degree $
L $ (the domain of the elliptic CAZP is the Riemann sphere), the
function (\ref{flat}) with probability one is an entire function,
and the function (\ref{hyperbolic}) with probability one is
analytic in the unit disc.

We introduce unified notation: $ \mathcal M $ for the domain of
the CAZP, and $ \Gamma $ for the symmetry group of $ \mathcal M$.
Then {\em CAZP is a unique $\Gamma$-stationary random zero
process}. Here, $\Gamma$-stationarity means that for any
$\gamma\in\Gamma$ the point processes $\mathcal Z_f$ and $\mathcal
Z_{f\circ \gamma}$ have the same distribution. Uniqueness means
that CAZP is the only $\Gamma$-stationary process on $\mathcal M$
among the random zero sets of Gaussian analytic functions.

Having explicit formulas \eqref{elliptic}, \eqref{flat},
\eqref{hyperbolic}, it is very easy to prove $\Gamma$-stationarity
and uniqueness. It suffices only to check that
\begin{equation}\label{eq*}
\mathbb E n_{f_L} = L \cdot m^*\,,
\end{equation}
where $m^*$ is a normalized $\Gamma$-invariant area measure on
$\mathcal M$. Then, by Calabi's rigidity, $\mathcal Z_{f_L}$ is
$\Gamma$-stationary, and, again by Calabi's rigidity, $f_L$ is
unique. Verification of \eqref{eq*} is a straightforward
application of the Edelman-Kostlan formula. For example, in the
flat case, $\|f_L\|^2(z) = \exp (L |z|^2)$, and
\[
n_{f_L} = \frac1{2\pi} \Delta (L|z|^2/2) = L \cdot \frac1{\pi}
m\,.
\]
We see that in the flat case the normalized area $m^*$ equals
$\frac1{\pi} m$. It is also worth mentioning that canonical
isometric embeddings of $\mathcal M$ into projective Hilbert
spaces corresponding to the Gaussian analytic functions
\eqref{elliptic}, \eqref{flat} and \eqref{hyperbolic} are well
known in geometry and physics.

By \eqref{eq*}, the parameter $L$ equals the mean number of random
zeroes per unit area on $\mathcal M$. In what follows, we shall
consider the asymptotic behaviour of the random zero processes
$\mathcal Z_{f_L}$ in the `{\em large intensity limit}'
$L\to\infty$. Some features for all three canonical models are
similar, some are different. Due to compactness, the elliptic
model sometimes is simpler to analyze. On the other hand, the
hyperbolic model has additional intriguing features.

To fix ideas, we shall concentrate on the flat model. In this
case, introducing $L$, we just make a homothety of the plane with
coefficient $r=\sqrt L$. This makes the flat case more
transparent\footnote{The scaling $z=z_0 + \frac{w}{\sqrt{L}}$
flattens out the elliptic and hyperbolic geometry as $L\to\infty$.
In this limit, the entire function \eqref{flat} is a locally
uniform limit of the functions \eqref{elliptic} and
\eqref{hyperbolic}, and the flat CAZP appears as a scaling limit
of the other two CAZP models. This is the motivation for an
advanced theory developed by Bleher, Shiffman and Zelditch in
\cite{BSZ}.}. Thus, we do not need parameter $L$ anymore, and we
consider the asymptotic zero distribution of the Gaussian entire
function of order two
\[
f(z) = \sum_{k=0}^\infty \zeta_k \frac{z^k}{\sqrt{k!}}\,.
\]
The function $f(z)$ can be viewed as a Gaussian counterpart of the
Weierstrass $\sigma$-function
\[
\sigma (z) = z \prod_{\omega\in{\sqrt\pi} \mathbb Z^2\setminus
\{0\}} \left( 1 - \frac{z}{\omega} \right) e^{z/\omega + 1/2
(z/\omega)^2}\,.
\]
Indeed, the random function $|f(z)|e^{-|z|^2/2}$ has a stationary
distribution, while the function $|\sigma (z)|e^{- |z|^2/2}$ has
periods $\sqrt\pi$ and $i \sqrt\pi$.

\section{Linear statistics}\label{sect3} Given a test-function $h\colon
\mathbb C \to \mathbb R$ with a compact support, consider the
random variable
\[
Z_r(h) = \int h(\tfrac{z}{r})\, dn (z)\,, \qquad \mathbb E Z_r (h)
= \frac{r^2}{\pi} \int h\, dm\,,
\]
$n$ is a counting measure of the flat CAZP with intensity $L=1$,
$m$ is the area measure. We are interested in the asymptotic
behaviour of $ Z_r(h) $ when $r\to\infty$. The size of
fluctuations of $Z_r(h)$ depends on the smoothness of the
test-function $h$.

\subsection{Smooth linear statistics}
\begin{thm}[\cite{ST}]\label{thm.var} Let $h$ be a $C^2$-function on
$\mathbb C$ with a compact support. Then
\begin{equation}\label{eq.var}
\operatorname{Var} Z_r(h) = \frac{\kappa + o(1)}{r^2} \, \|\Delta
h \|^2_2\,, \qquad r\to\infty\,,
\end{equation}
where $\kappa$ is a positive numerical constant. The random
variables
\[
\frac{r}{\sqrt{\kappa} \|\Delta h\|_2}\, \left( Z_r(h) -
\frac{r^2}{\pi} \int h\, dm \right)
\]
converge in distribution to the standard Gaussian law $\mathcal
N(0;1)$ for $r\to\infty$.
\end{thm}

Asymptotic formula \eqref{eq.var} first appeared in Forrester and
Honner~\cite{FH}. It is worth mentioning that
Theorem~\ref{thm.var} persists for the other two CAZP models in
the large intensity limit $L\to\infty$ \cite[Part I]{ST}.

It is instructive to compare \eqref{eq.var} with the size of
variations for a simple point process in the plane given by i.i.d.
Gaussian perturbations of the lattice. Consider the point process
\[
\mathcal S = \left\{\sqrt{\pi} (k+il) + \eta_{k,l}\colon
(k,l)\in\mathbb Z^2 \right\}
\]
where $\eta_{k,l}$ are independent standard complex Gaussian
random variables. In this case, $\operatorname{Var} S_r(h) \sim
\operatorname{const} \|\nabla h\|_2^2$, for $r\to\infty$. This is
rather different from \eqref{eq.var}. Asymptotic similarity to the
flat CAZP $\mathcal Z_f$ can be achieved by inventing special
correlations between the perturbations $\eta_{k,l}$.
\footnote{Lattice points are aggregated into clusters and each
cluster scatters in a special (equiangular and equidistant) way
\cite[Part~I, Introduction]{ST}.} In Section~\ref{sect4}, we
return to the idea of the flat CAZP as a perturbed lattice.

The proof of Theorem~\ref{thm.var} starts with Green formula
\[
Z_r(h) - \frac{r^2}{\pi} \int h\, dm  = \frac1{2\pi}\, \int \log
|f_r^*| \, \Delta h\, dm\,, \qquad f_r^*(z) =
\frac{f(rz)}{\sqrt{\text{Var} f(rz)}}\,.
\]
The RHS is a non-linear functional on a Gaussian process $f_r^*$.
The rest is based on the method of moments \'a la Breuer and Major
\cite{BM}: we expand the function $\zeta\mapsto \log|\zeta|$ in
Hermite polynomials in the space $L^2_{\mathbb C}(e^{-|\zeta|^2})$
(the Wick expansion), and evaluate the moments of $Z_r(h)$ using
the combinatorial diagram technique.

\subsection{Number of chaotic analytic zero points} Let
$\Omega\subset \mathbb C$ be a bounded domain with a piecewise
smooth boundary. We are interested in the asymptotic behaviour of
the random variable $n(r\Omega) = Z_r(\done_\Omega)$ for
$r\to\infty$. Forrester and Honner \cite{FH} argued that the
technique developed by Martin and Yal\c{c}in \cite{MY} for
studying the Gibbs states of infinite systems of charged particles
applied to the flat CAZP gives
\[
\text{Var}\, n(r\Omega) = r\cdot(\tau +
o(1))\,\text{Length}(\partial\Omega)\,, \qquad r\to\infty\,,
\]
$\tau$ is a positive numerical constant. This is consistent with
the idea that the variation of the number of points in $r\Omega$
should behave like the number of points in the `strip' of constant size
around the boundary $\partial (r\Omega)$.

Presumably, the method of Martin and Yal\c{c}in also yields that
the random variables
\[
\frac{n(r\Omega) - \pi^{-1} r^2 m(\Omega)}{\sqrt{r\cdot \tau \,
\text{Length}(\partial\Omega)}}
\]
converge in distribution to $\mathcal N(0;1)$ for $r\to\infty$.

It would be interesting to find a counterpart of the law of the
iterated logarithm; i.e. to find a function $\phi (r)$ such that
{\em with probability one}
\[
\limsup_{r\to\infty} \frac{|n(r) - r^2|}{\sqrt{r} \phi (r)} = 1\,.
\]
Here $n(r)=n\left(\{|z|\le r\}\right)$.

\subsection{The `hole probability' and large deviations}
The next theorem proves an estimate conjectured by Yuval Peres:
\begin{thm}[\cite{ST}]\label{thm.hole} For $r\ge 1$,
\begin{equation}\label{eq.60}
e^{-c_1 r^4} \le  \mathbb P\left( n(r) = 0 \right) \le e^{-c_2
r^4}\,,
\end{equation}
where $n(r)=n(\{|z|\le r\})$, and  $c_1$ and $c_2$ are positive
numerical constants.
\end{thm}

It would be interesting to check whether there exists the limit
\[
\lim_{r\to\infty} \frac{\log^- \mathbb P \left( n(r)=0
\right)}{r^4}\,,
\]
and (if it does) to compute its value.

The lower bound in \eqref{eq.60} is obtained by an explicit
construction. The upper bound follows from
\begin{thm}[\cite{ST}]\label{thm.large}
For any $\delta\in (0, \frac14]$ and $r\ge 1$,
\[
\mathbb P \left( \left| \frac{n(r)}{r^2} -1 \right| \ge \delta
\right) \le e^{-c(\delta) r^4}\,.
\]
\end{thm}

The proof of Theorem~\ref{thm.large} uses tools from the entire
functions theory. First, we show that with very high probability
$\log \max_{r\mathbb D} |f|$ is close to $r^2/2$. Then, estimating
$\log |f|$ from below, we show that with very high probability the
average
\[
\int_0^{2\pi} \log |f(re^{i\theta})|\, \frac{d\theta}{2\pi}
\]
is also close to $r^2/2$. From this, using Jensen's formula, we
deduce Theorem~\ref{thm.large}.

Theorems~\ref{thm.hole} and \ref{thm.large} are consistent with
the results known for a one component Coulomb system of charged
particles  of one sign embedded into a uniform background of the
opposite sign, Jancovici, Lebowitz and Magnificat \cite{JLM}. It
would be good to understand the asymptotic behaviour of the random
variable $r^{-\alpha}\left[ n(r) - r^2 \right]$ for $r\to\infty$
and $\alpha\ge 1/2$. At present, we understand the extreme case
$\alpha = \frac12$, and (not completely) the case $\alpha \ge 2$.
A plausible guess (motivated by \cite{JLM}) is
\[
\lim_{r\to\infty} \frac{\log\log \mathbb P \left( |n(r) - r^2| \ge
r^\alpha \right)}{\log r} =
\begin{cases}
2\alpha -1, &\frac12 \le \alpha \le 1; \\
3\alpha -2, &1\le\alpha\le 2;\\
2\alpha, &\alpha \ge 2\,.
\end{cases}
\]
In the first case, the normalized charge $|n(r) - \frac{r^2}2|$
grows slower than the perimeter of the disc, in the second case it
grows faster than the perimeter but slower than the area, in the
last case (so called `overcrowding') it grows faster than the
area. According to the philosophy of \cite{JLM}, this should lead
to a change of the asymptotic regime at $\alpha =1$ and $\alpha
=2$. The technique we developed for the proof of
Theorem~\ref{thm.hole} helps to analyze the case $\alpha>2$. The
other two cases seem to require a different technique.

\section{Flat chaotic analytic zero points as a perturbed
lattice} \label{sect4}

How evenly do the flat CAZP spread over the plane? In the
euclidean case we have a very natural system of points spread
evenly throughout the plane, namely, the lattice points. Instead
of comparing the random counting measure $n_f$ with its average
$\frac1{\pi} m$, we consider the flat CAZP as a perturbed lattice
$\{\sqrt{\pi}(k+il)+\xi_{k,l}\colon k,l\in \mathbb Z \}$, for some
(dependent) complex-valued random variables $\xi_{k,l}$. We may
hope for fast decay of the tail probabilities $\mathbb P\left(
|\xi_{k,l}|\ge \lambda \right)$ for large $\lambda$, uniformly in
$(k,l)\in\mathbb Z^2$. The uniformity becomes trivial if the
distribution of $(\xi_{k,l})$ is invariant under lattice shifts.
We treat the random variables $\xi_{k,l}$ as measurable functions
on the space $\Omega = \mathbb C^{\mathbb Z^2}$ of two-dimensional
arrays $\xi\colon \mathbb Z^2\to\mathbb C$.
\begin{thm}[\cite{ST}]\label{thm.lattice}
There exists a probability measure $\mathbb P$ on $\Omega$ such
that
\newline\smallskip\noindent
(i) $\mathbb P$ is invariant under the shifts of $\mathbb Z^2$;
\newline\smallskip\noindent
(ii) the random set $\mathcal S = \left\{\sqrt{\pi}(k+il) +
\xi_{k,l}\colon (k,l)\in\mathbb Z^2 \right\}$ is distributed like
the flat CAZP $\mathcal Z$;
\newline\smallskip\noindent
(iii) $\mathbb E\left( e^{\epsilon |\xi_{0,0}|^2}\right)<\infty $
for some $\epsilon>0$.
\end{thm}

This result gives no information about correlation between
$\xi_{k,l}$. Probably, they can be chosen to be nearly independent
on large distances.

The proof of Theorem~\ref{thm.lattice} does not use the Gaussian
nature of the flat CAZP but only the uniform boundedness of the
exponential moment of the `random potential' $u(z)=\log |f(z)| -
\frac12 |z|^2$. The main ingredients of the proof are the
M.~Hall's `marriage lemma' (needed to match the flat CAZP with the
lattice $\sqrt{\pi}\mathbb Z^2$), and a potential theory lemma
which can be useful in other discrepancy problems:
\begin{lemma}\label{lemma.transp}
Let $u$ be a bounded delta-subharmonic function on $\mathbb C$ (i.e.,
a difference of two subharmonic functions), and let $\Delta u = \mu
- m$ in the distributional sense, $\mu$ is a non-negative measure.
Then for any bounded Borel set $E\subset \mathbb C$
\[
\mu(E) \le m(E_{+t}) \qquad \text{and} \qquad m(E) \le \mu (
E_{+t})\,,
\]
where $E_{+t} = \{z\colon \operatorname{dist}(z, E)\le t\}$ is a
$t$-vicinity of $E$,
$t=\operatorname{const}\,\|u\|_{\infty}^{1/2}$.
\end{lemma}

The boundedness of $u$ is too strong for applications. It can be
easily weakened by convolving $u$ with a smooth convolutor
supported by an appropriate disc.

The two ingredients described above alone help to prove only a
local result in the spirit of Theorem~\ref{thm.lattice}.
Globalization is still a problem: after smoothing, the random
potential $u(z)=\log |f(z)| - \frac12 |z|^2$ is almost surely
unbounded. Rare fluctuations appear somewhere on the infinite
plane $\mathbb C$, though probably far from the origin. To achieve
some locality, we introduce a special adaptive metric $\rho$ on
$\mathbb C$ using H\"ormander's construction
\cite[Section~1.4]{Hor}. This metric $\rho$ is small where the
potential is large. Then we use a version of
Lemma~\ref{lemma.transp} making use of $\rho$-neighbourhoods
instead of euclidean ones.

It is worth mentioning that in \cite{AKT} Ajtai, Koml\'os and
Tusn\'ady studied high probability matchings of a system of $N^2$
independent random points $\Lambda = \{\lambda_1, ...\,
\lambda_{N^2}\}$ uniformly distributed in the square $[0, N]^2
\subset \mathbb R^2$ with the grid $\{\omega_1, ...\,
\omega_{N^2}\} = [0, N)\cap \mathbb Z^2$. They considered the
average transportation
\[
T(\Lambda) := \min_{\pi} \frac1{N^2} \sum_{1\le i \le N^2}
|\lambda_i - \omega_{\pi(i)}|\,,
\]
the minimum is taken over all permutations $\pi$ of $\{1, 2, ...\,
N^2\}$. Then with high probability
\begin{equation}\label{eq.AKT}
\text{const} \sqrt{\log N} \le T(\Lambda) \le \text{Const}
\sqrt{\log N}\,.
\end{equation}
For related results see Leighton and Shor~\cite{LS}, and
Talagrand~\cite{T1}. In our global set-up we deal with {\em
infinite} measures in the plane. Then, according to the lower
bound in \eqref{eq.AKT}, for {\em any} matching the average
transportation distance tends to infinity in the $N\to\infty$
limit. This leaves no hope for the finite average distance
matching between the Poissonian point process in $\mathbb R^2$ and
a lattice, even without a quantitative estimate (iii).

It would be interesting to find a hyperbolic counterpart of
Theorem~\ref{thm.lattice}.

\section{Hyperbolic CAZP}\label{sect5}

The hyperbolic CAZP, in contrast to the other models, depends on
two intrinsic parameters: the mean density of zeroes per unit
hyperbolic area, and the size of the sets on which we count the
number of points. For instance, $n_L(D_r)$ is a number of the
hyperbolic CAZP with intensity $L$ in the hyperbolic disc of
radius $r$. This leads to different asymptotic regimes, and makes
the hyperbolic model richer than the other two. It is also worth
mentioning that, for different values of the intensity $L$, the
Gaussian analytic function \eqref{hyperbolic} exhibits completely
different patterns in the asymptotic behaviour when $z$ approaches
the boundary of the unit disc Kahane \cite[Chapter~13]{Kahane}.

An interesting observation  by Diaconis and Evans \cite{DE}, and
Peres and Vir\'ag \cite{PV} says that the real part of the
hyperbolic random function \eqref{hyperbolic} up to a constant
term is a Poisson integral of a Gaussian random noise on the unit
circle. In the case $L=1$ this is classical white noise
\cite[p.11]{PV}, the case $L=2$ corresponds to the Gaussian
process on the Dirichlet space $H_2^{1/2}(\mathbb T)$
\cite[Example~5.6]{DE}.

\subsection{An exactly solvable model} Here, we discuss a recent
finding of Peres and Vir\'ag \cite{PV} which pertains to the case
$L=1$. Recall that the $k$-point correlation function
$p(z_1,\,...\,z_k)$ of a random point process is
\[
p(z_1,\,...\,z_k) = \lim_{\epsilon\to 0} \frac{p_\epsilon
(z_1,\,...\,z_k)}{(\pi \epsilon^2)^k}\,,
\]
where $p_\epsilon(z_1,\,...\,z_k)$ is the probability that each
disc $\{|z-z_j|\le \epsilon\}$, $1\le j\le k$, contains at least
one point of the process. For the random zero processes of a
Gaussian analytic function, the limit on the RHS always exists. An
equivalent definition says that $p(z_1, \, ...\, z_k)$ is a
density of the $k$-point correlation measure \eqref{eq.mu} with
respect to the Lebesgue measure. The correlation measure also
contains singular terms supported by the large diagonal; these
terms are expressed via $j$-point correlation functions with
$j<k$. Thus, the random zero process can be described by its
correlation functions.

Using Hannay's formulas \cite{Hannay96}, Peres and Vir\'ag proved
\begin{thm}[Peres-Vir\'ag \cite{PV}] The correlation function of
the hyperbolic CAZP with intensity $L=1$ is
\[
p(z_1,\, ...\,z_k) = \pi^{-n} \det\left[
\frac1{(1-z_i\bar{z_j})^2}\right]_{1\le i,j \le k}\,.
\]
\end{thm}

This remarkable identity  makes  the hyperbolic CAZP with  $L=1$
an `exactly solvable model' among all CAZP. \footnote{Peres and
Vir\'ag observe that this is the only {\em determinantal} process
among CAZP.} In particular, it yields amazingly simple explicit
expressions for the distribution of the number of zeroes $n(\rho)$
in the disc $\{|z|\le\rho\}$ and for the asymptotics of the `hole
probability':
\begin{cor}[\cite{PV}]
Let $\mathcal Z$ be the hyperbolic CAZP with intensity $L=1$. Then
\newline\noindent (i) $n(\rho)$ has the same distribution as
$\sum_{j=1}^\infty X_j$ where $\{X_j\}$ is a sequence of
independent Bernoulli random variables with $\mathbb P(X_j=1) =
\rho^{2j}$;
\newline\noindent (ii) for $\rho\to 1$
\begin{equation}\label{eq.70}
\mathbb P \left( n(\rho) = 0 \right) = \exp \left[ - \frac{\pi^2 +
o(1)}{1-\rho} \right];
\end{equation}
\newline\noindent (iii) the ratio
\[\frac{n(\rho) - \mathbb E n(\rho)}{\sqrt{\operatorname{Var}
n(\rho)}}
\]
(with $\mathbb E n(\rho) = \frac{\rho^2}{1-\rho^2}$
and $\operatorname{Var} n(\rho) = \frac{\rho^2}{1-\rho^4}$)
converges in distribution to the standard Gaussian law $\mathcal
N(0;1)$ for $\rho\to 1$.
\end{cor}

In this case, $\operatorname{Var} n(\rho)$ has the same order of
magnitude as $\mathbb E n(\rho)$ whilst in the flat case the
variance grows only as a square root of the mean. This naturally
reflects the difference between the hyperbolic and euclidean
geometries: in hyperbolic geometry the perimeter grows like the
area, and much more random zeroes are located near the boundary
circumference.

We are not aware of counterparts of (ii) and (iii) for the
hyperbolic CAZP with $L\ne 1$.

\section*{Loose ends} In this lecture we've only touched the
`ground level' of the theory. There are plenty of interesting and
deep developments. Among them are
\begin{itemize}
\item[$\bullet$] scaling limits of
zeroes of random polynomials (Ibragimov and Zeitouni \cite{IZ} and
Shiffman and Zelditch~\cite{SZ0} ) and of random holomorphic
sections of high powers of Hermitian line bundles (Bleher,
Shiffman and Zelditch~\cite{BSZ});
\item[$\bullet$] solutions of random systems
of algebraic equations, including sparse systems (see Edelman and
Kostlan~\cite{EK}, Malajovich and Rojas~\cite{MR}, Shiffman and
Zelditch~\cite{SZ} and references therein);
\item[$\bullet$]
distribution of real zeroes of random real polynomials
(Maslova~\cite{Maslova}, Dembo, Poonen, Shao and
Zeitouni~\cite{DPSZ}, Bleher and Di~\cite{BD}, Aldous and
Fyodorov~\cite{AF});
\item[$\bullet$] links with the zero distribution
of chaotic eigenfunctions (Nonnenmacher and Voros~\cite{NV},
Hannay~\cite{Hannay96}) and with the distribution of eigenvalues
of random matrices with independent complex Gaussian entries
(Forrester and Honner~\cite{FH}, Diaconis and Evans~\cite{DE},
Dennis and Hannay~\cite{DH} and references therein),
\end{itemize}
and each of them deserves a special lecture (cf. \cite{Zelditch}).
But these are different stories to be told by other people.

\paragraph*{Acknowledgement} For the last four years, I have been
enjoying collaboration with Boris Tsirelson on the subject of this
lecture. I am grateful to him for his encouragement and patience.
I am also grateful to F\"edor Nazarov, Leonid Pastur, Yuval Peres,
Leonid Polterovich, Zeev Rudnick, Bernard Shiffman, Peter Yuditskii,
and Steve Zelditch for useful discussions.


\end{document}